\documentclass[a4paper,12pt]{article}
\usepackage[utf8]{inputenc}
\usepackage[russian]{babel}
\usepackage{amsmath,amsthm,amssymb,amscd,array}
\usepackage{hyperref}
\usepackage{titlesec}

\renewcommand{\geq}{\geqslant}
\renewcommand{\leq}{\leqslant}
\renewcommand{\epsilon}{\varepsilon}

\newcommand{\w}{\ensuremath{\mathop{W}\limits^\circ {}}}

\theoremstyle{definition}
\newtheorem*{defn}{Определение}

\titleformat{\section}
[display]
{\centering\normalfont\large\bfseries}
{}
{0pt}
{\centering\large\bfseries}

\titlespacing*{\section} {10pt}{-10pt}{10pt}

\makeatletter
\g@addto@macro\normalsize{%
\setlength{\abovedisplayskip}{2pt}
\setlength{\belowdisplayskip}{2pt}
\setlength{\abovedisplayshortskip}{1pt}
\setlength{\belowdisplayshortskip}{1pt}
}
\makeatother

\begin{document}

	\title{\bf О НЕПЕРЕСЕЧЕНИИ СПЕКТРОВ МИНИМИЗИРУЕМЫХ ФУНКЦИОНАЛОВ В ПРОСТРАНСТВАХ $\w_2^n$, $\w_2^{n+1}$, $\w_2^{n+2}$}
		\date{}
	\author{Андрей Михайлович Минарский}
	
	\maketitle
	
\vspace{-2em}
\begin{center}
Физико-Техническая Школа, СПб	\\
Сондужская Высшая Школа, Тотьма
\end{center}

\begin{center}
\bf Аннотация
\end{center}

\begin{small}
Рассмотрены спектры собственных значений функционала
$$\Phi(u)=\frac{\left\langle u^{(n)}u^{(n)}\right\rangle}{\left\langle u^{(n-p)}u^{(n-p)}\right\rangle}$$
в $\w_2^n$ при разных $n$. Показано, что спектры на четных функциях в $\w_2^n$ и $\w_2^m$ не пересекаются при $m=n+1, n+2$. Написаны необходимые условия возможности пересечения при $\Delta=m-n>2$.
\end{small}
		
\bigskip


\section{Введение}

Пусть дано (соболевское) пространство вещественных функций $$\w_2^n=\left\{u\in C^{n-1}\,[0, 1] \,\,|\, u^{(j)}(0)=u^{(j)}(1)=0,\;j=0,\ldots, n-1,\;u^{(n)}\in L_2\,(0, 1)\right\}.$$
Введем стандартную норму $$\|u\|_n^2 = \left\langle u^{(n)}u^{(n)}\right\rangle \text{ где } \,\,\langle a\,b\rangle=\int\limits_0^1 a(t)b(t)\,dt.$$

Рассмотрим функционал 
\begin{align}
	\label{eq1}
	\Phi(u) = \frac{\left\langle u^{(n)}\;u^{(n)}\right\rangle}{\left\langle u^{(n-p)}\;u^{(n-p)}\right\rangle}.
\end{align}

Задача нахождения минимума
\begin{align*}
	\Lambda_n=\min_{u\in\w_2^n}\Phi(u)
\end{align*}
при варьировании с $\delta u=\epsilon v$ приводит к уравнению
\begin{align}
	\label{eq211}
	\left.\Phi'_\epsilon\right|_{\epsilon=0}=0=\left\langle u^{(n)}v^{(n)}\right\rangle - \Lambda\left\langle u^{(n-p)}v^{(n-p)}\right\rangle\quad \text{для}\;\forall v\in\w_2^n.
\end{align}

Обозначим $$\Omega^n=C^r\,[0, 1]\mathop{\cap}\w_2^n,$$ где $r\geq2n$ --- наибольшая используемая в дальнейшем исследовании производная.

При рассмотрении $\Phi(y)$ на $y\in\Omega^n$, из \eqref{eq211} получим
\begin{align}
	\label{eq3}
	\left\langle v L^{2n} y\right\rangle=0 \text{ для всякого }  v\in\Omega^n,
\end{align}
где оператор
\begin{align}
	\label{eq2}
	L^{2n}=L^{2n}(\Lambda)=(-1)^n d^{2n} - \Lambda(-1)^{n-p}d^{2n-2p},\quad d=\frac{d}{dt}.
\end{align}

Если функция $z$ есть решение \eqref{eq3}, то она имеет вид
\begin{align}
	\label{eq4}
	z=R+P_{2n-2p-1},
\end{align}
где $R$~--- часть, удовлетворяющая $\left[d^{2p}-(-1)^p\Lambda\right]R=0$, а $P$~--- полином. В дальнейшем будем писать $P_k$ для полиномов с $\deg P_k\leq k$.
\\

Обозначим через $z_n$ собственную функцию (с.ф.) $L^{2n}$: $L^{2n}z_n=0$, а через $\overline{\mathstrut \Lambda}_n$~--- соответствующее собственное значение (с.з.): $L^{2n}\left(\overline{\mathstrut\Lambda}_n\right)z_n=0$.  Пусть $S_n=\left\{\overline{\mathstrut \Lambda}_n\right\}$~--- спектр с.з., а $\Lambda_n$~--- минимальное с.з.

Спектр с.з. для оператора $L^{2n}$ в $\Omega^n$ совпадает со спектром с.з. $\Phi(u)$ в $\w_2^n$. Далее будем изучать ситуацию в $\Omega^n$. Замена $t=\frac{x+1}2$ превращает $[0, 1]\rightarrow[-1, 1]$, при соответствующем переобозначении $d$ в \eqref{eq2} и изменении с.з.
\\

С.ф. разделим на симметричные $z_{n, s}=z_n(x)+z_n(-x)$ и антисимметричные $z_{n, a}=z_n(x)-z_n(-x)$. Соответствующие с.з. и спектры обозначим $\overline{\mathstrut \mathstrut \Lambda}_{n, s}$, $\overline{\mathstrut \mathstrut \Lambda}_{n, a}$; $S_{n, s}$, $S_{n, a}$.
 
Очевидно, что $S_{n, a}=S_{n+1, s}$.\footnote{Следует из того, что $z_{n, a}^{(-1)}=\int_{-1}^x z_{n, a}(\xi)\,d\xi$ есть некоторая $z_{n+1, s}$ и обратно $\forall z_{n+1, s}^{(1)}$ в силу (\ref{eq4}, \ref{eq5}) есть $z_{n, a}$.} В силу $\Omega^{n+1}\mathop{\subset}\Omega^n$ имеем $\Lambda_{n, s}\leq\Lambda_{n+1, s}=\Lambda_{n, a}$. В дальнейшем, не оговаривая, рассматриваем только симметричные функции.

Заметим, что
\begin{align}
	\label{eq5}
	z_{ns}=R+P_{2n-2p-2}.
\end{align}

Ниже везде в обозначении $L^{2k}z_l$ подразумевается, что $\Lambda$ в $L^{2k}$ есть с.з. оперируемой функции $z_l$ для $L^{2l}$: $L^{2l}(\Lambda)z_l=0$.


\section{Камень в оперируемом}

В разделе доказано $S_{n-1, s}\mathop{\cap}S_{n, s}=\varnothing$.

\begin{defn}
Пусть $n>p$. Назовем
\begin{align}
	\label{eq6}
	c_n=L^{2n-2}z_n
\end{align}
\emph{камнем} в $z_n$.
\end{defn}
	
\noindent\textbf{Лемма (о <<неизымаемости>> камня):}
\textit{Если $c_n=0$, то $z_n=0$.}

\noindent\textbf{Доказательство.}
Имеем:
\begin{align}
	\label{eq7}
	\left\langle z_n L^{2n}z_n\right\rangle - \lambda^2\left\langle z_n L^{2n-2}z_n\right\rangle = -\lambda^2 c_n \left\langle z_n\right\rangle,
\end{align}
где $\langle a\rangle=\langle 1, a\rangle$.

Обозначим $\sigma = i\cdot d = i\cdot \frac{d}{dx}$ (<<импульс>>). Тогда $$L^{2n}-\lambda^2 L^{2n-2}=\left(\sigma^2-\lambda^2\right)L^{2n-2}.$$
	
Пусть
\begin{align}
	\label{eq71}
	\sigma^{2p}-\overline{\mathstrut \Lambda_n} = \left(\sigma^2-\lambda_0^2\right)\cdot\prod\limits_{i=1}^{p-1}\left(\sigma-\lambda_i\right)\cdot\left(\sigma-\overline{\mathstrut \lambda_i}\right),
\end{align}
где $\lambda_0$~--- вещественный, $\lambda_i$ и $\overline{\mathstrut \lambda_i}$~--- сопряженные корни $\sigma^{2p}-\overline{\mathstrut\Lambda}_n$ как полинома от $\sigma$.

Выберем $\lambda=\lambda_0$. Из \eqref{eq2} и \eqref{eq71} очевидно:
\begin{align}
	\label{eq8}
	L^{2n}-\lambda^2 L^{2n-2}=\overline{\mathstrut A}(\sigma)A(\sigma),\quad \text{где}\;A=\sigma^{n-p-1}\left(\sigma^2-\lambda^2\right)\cdot\prod\limits_{i=1}^{p-1}\left(\sigma-\lambda_i\right),
\end{align}
а $\overline{\mathstrut A}$, как полином от $\sigma$, комплексно сопряжен $A$.
		
Из (\ref{eq7}, \ref{eq8}) получаем
\begin{align}
	\label{eq10}
	\left\langle z_n\overline{\mathstrut A}(\sigma)A(\sigma)z_n\right\rangle = -\lambda^2 c_n \left\langle z_n\right\rangle.
\end{align}
		
Отсюда, в силу эрмитовости в $\Omega^n$ операторов $\sigma^k$ при $k\leq2n$, если $c_n=0$:
\begin{align}
	\label{eq111}
	A(\sigma)z_n=\overline{\mathstrut A}(\sigma)z_n=0.
\end{align}
		
Тогда $R$ (см. \eqref{eq5}) содержит лишь $\cos\lambda x$ и $\left(\sigma^2-\lambda^2\right)z_n = Q$~--- полином. В силу (\ref{eq8}, \ref{eq111}) $\deg Q\leq n-p-2\leq n-3$, однако $Q\in\Omega^{n-2}\Rightarrow Q=0.$ Остающаяся в $z_n$ часть $\cos\lambda x$ исчезает в $\Omega^n$ при $\forall n\geq2$.\hfill$\square$
	
\noindent\textbf{Докажем} $\overline{\mathstrut \Lambda}_{n-1}\neq\overline{\mathstrut \Lambda}_n.$ Имеем:
\begin{align}
	\label{eq11}
	\left\langle z_n^{(n-1)} z_{n-1}^{(n-1)}\right\rangle =\left[
	\begin{array}{l}
		\overline{\mathstrut \Lambda}_{n-1}\langle\cdots\rangle + \left\langle z_n L^{2n-2} z_{n-1}\right\rangle = \overline{\mathstrut \Lambda}_{n-1}\langle\cdots\rangle\\
		\overline{\mathstrut \Lambda}_n\langle\cdots\rangle + \left\langle z_{n-1} L^{2n-2} z_n\right\rangle = \overline{\mathstrut \Lambda}_n\langle\cdots\rangle + c_n\left\langle z_{n-1}\right\rangle
	\end{array}
	\right.,
\end{align}
где $\langle\ldots\rangle = \left\langle z_n^{(n-p-1)} z_{n-1}^{(n-p-1)}\right\rangle$.
		
При $\overline{\mathstrut \Lambda} _n= \overline{\mathstrut \Lambda}_{n-1}$ имеем $c_n\left\langle z_{n-1}\right\rangle = 0$. В силу (\ref{eq7}, \ref{eq10}) и леммы о камне либо $z_n$, либо $z_{n-1}$ исчезают.\footnote{Условие $n>p$ обеспечивается существованием $z_{n-1}$. Не попавший прямо под доказательство случай $n-1=p$, $\left\langle z_{n-1}\right\rangle=0$ устраняется тем, что тогда $z_{n, a}=z_{n-1, s}^{(-1)}=\int_{-1}^x z_{n-1, s}(\xi)\,d\xi$ и $z_{n+1, s}''=z_{n, a}'=z_{n-1, s}$, откуда $c_{n+1}=L^{2n}z_{n+1, s}=L^{2n-2}z_{n-1, s}=0$.} 
	
Итак, доказано $S_{n-1, s}\mathop{\cap}S_{n, s}=\varnothing$ при $\forall n$. В частности, отсюда $\Lambda_{n-1}<\Lambda_n$.

\noindent\textbf{Утверждение (о простоте спектра в $\Omega^n$):} Если у $L^{2n}$ есть две линейно независимые с.ф. в $\Omega^n$: $z_{1, n}$ и $z_{2, n}$, то
\begin{align}
	\Lambda_{1, n}\ne\Lambda_{2, n}.
\end{align}

Действительно, пусть не так. Тогда их линейная комбинация $z=\alpha\cdot z_{1, n}+\beta\cdot z_{2, n}$, где $\alpha=z_{2, n}(1)$, $\beta=-z_{1, n}(1)$, такова, что $z\in\Omega^{n+1}$, $L^{2n}z=0\Rightarrow z=0$ в силу леммы о камне.
	

\section{Соотношения между камнями и моментами}

Далее обозначаем $\Lambda_n$ вместо $\overline{\mathstrut\Lambda}_n$. В разделе доказано $\Lambda_n\neq\Lambda_{n+\Delta}$ при $\Delta = 1, 2$. По\-лу\-че\-ны некоторые требования для возможности выполнения равенства при б\' oльших $\Delta$.

Назовем (неполным) \emph{полиномом камней}:
\begin{align}
	h_n^k&=(-1)^k L^{2n-2k-2}z_n; \notag \\
	\label{eq12}
	h_n^k&=c_{0n}\frac{z^{2k}}{(2k)!}+c_{1n}\frac{z^{2k-2}}{(2k-2)!}+\ldots+c_{kn}; \\
	h_n^k&=0,\quad k<0. \notag
\end{align}

Очевидно $d^2 h_n^k=h_n^{k-1}$.

Пусть $m>n$; обозначим $\mu=n+2\left[\frac{m-n}2\right]$. Выполняя для $\left\langle z_n^{(n-k)}z_m^{(\mu-k)}\right\rangle$ то же, что в \eqref{eq11}, с учетом \eqref{eq12} и сдвигая $k$ на $-1-[\Delta/2]$, получим
\begin{align}
	\label{eq13}
	\left\langle z_m h_n^k\right\rangle-(-1)^\Delta\left\langle z_n h_m^{k+\Delta}\right\rangle = \left(\Lambda_m - \Lambda_n\right)\cdot(-1)^k\left\langle\ldots\right\rangle
\end{align}
для $\forall k\geq -1-\left[\Delta/2\right]$, $\Delta=m-n$, где $\left\langle\ldots\right\rangle=\left\langle z_n^{(n-p-k-1)}z_m^{(n-p-k-1)}\right\rangle$.

Обобщая \eqref{eq8}:
\begin{align}
	\left(\sigma^2-\lambda^2\right)L^{(2n-2k-2)}=\overline{A_{n-k}}A_{n-k}, \notag\\
	\text{ где }\;A_{n-k}=\sigma^{n-p-k-1}\left(\sigma^2-\lambda^2\right)\cdot\prod\limits_{i=1}^{p-1}\left(\sigma-\lambda_i\right).
\end{align}

Тогда
\begin{align}
	\label{eq15}
		\left\langle z_n\overline{A_{n-k}}A_{n-k}z_m\right\rangle=(-1)^{k-1}\left\langle z_n h_n^{k-1}\right\rangle - (-1)^k\lambda_n^2\left\langle z_n h_n^k\right\rangle>0, \,\,	\forall k\geq 0	&  \notag \\
		\text{и } (-1)^{k-1}\left\langle z_m h_m^{k-1}\right\rangle - (-1)^k\lambda_m^2\left\langle z_m h_m^k\right\rangle>0,  \text{ при } 	\lambda_n^2=\Lambda_n^{1/p},\;\lambda_m^2=\Lambda_m^{1/p}.
\end{align}

Для анализа условий \eqref{eq13} и \eqref{eq15} удобно представить их следующим образом: обозначим 
\begin{align}
	\label{eq16}
	\left\langle z_n \frac{x^{2k}}{(2k)!}\right\rangle=a_k,\quad \left\langle z_m \frac{x^{2k}}{(2k)!}\right\rangle=b_k,\quad c_{kn}=c_{k},\quad c_{km}=d_k,
\end{align}
где $a = \left(a_0, a_1, \dots\right)$, $b=\left(b_0, b_1, \dots\right)$~--- \emph{моменты} функций $z_n$ и $z_m$. Обозначим также
\begin{align}
	(f, g)_k=(-1)^k\left(f_k g_0+f_{k-1}g_1+\ldots+f_0 g_k\right);\quad(f, g)_k=0,\;k<0.
\end{align}

Тогда \eqref{eq13} и \eqref{eq15} перепишутся
\begin{align}
	\label{eq18}
	(b, c)_{k} - (a, d)_{k+\Delta} = \left(\Lambda_{n+\Delta}-\Lambda_n\right)\left\langle\ldots\right\rangle
	\intertext{для $\forall k: -1-\left[\frac{\Delta}2\right]\leq k\leq n-p-1$ и $\langle\ldots\rangle$ из \eqref{eq13}}
	\label{eq19}
	\left\{
	\begin{array}{ll}
		\epsilon_n(a, c)_{k-1}>(a, c)_k
		&
		\forall k\geq0
		\\
		\epsilon_m(b, d)_{k-1}>(b, d)_k
		&
		\epsilon_n=\Lambda_n^{-1/p},\;\epsilon_m=\Lambda_m^{-1/p}
	\end{array}
	\right.
\end{align}

Соотношения (\ref{eq18}, \ref{eq19}) закрывают возможность $\Lambda_n=\Lambda_{n+\Delta}$ при $\Delta=1, 2$ в силу $(ad)_0\cdot(bc)_0=(ac)_0\cdot(bd)_0\neq0$.

Неравенства \eqref{eq19} могут быть несколько усилены при $\Lambda_n=\Lambda_m(=\Lambda)$. Учитывая $$\left|\left\langle z_n\overline{\mathstrut A_{n-l}}A_{\mu-2k+l}z_m\right\rangle\right|^2<\left\langle z_n\overline{\mathstrut A_{n-l}}A_{n-l}z_n\right\rangle\cdot\left\langle z_m\overline{\mathstrut A_{\mu-2k+l}}A_{\mu-2k+l}z_m\right\rangle$$ и поступая, как при получении \eqref{eq15} и \eqref{eq19}, с учетом \eqref{eq18} получим
\begin{align}
	\label{eq21}
	(\alpha c)_l\cdot(\beta d)_{2k-l+\delta}>(\alpha d)_{k+q+\delta}\cdot(\beta c)_{k-q}
\end{align}
\begin{align*}
	\forall k:\;& 0\leq k\leq n-p-1+q,\\
	\forall l:\;& 0\leq l\leq n-p-1,\\
	          & -\delta\leq 2k-l\leq n-p-1+2q,
\end{align*}
где $q=[\Delta/2]$, $\delta=\Delta-2q$, $\Delta=m-n$. Обозначено $$\alpha=\left(\alpha_0, \alpha_1, \dots\right):\quad\alpha_0=a_0,\,\alpha_i=a_i+\epsilon a_{i-1}, i\in\mathbb N$$ и аналогично для $\beta$ и $b$, $\epsilon=\Lambda^{-1/p}$.

Ограничения сверху на $k$ и $l$ написаны на случай рассмотрения \emph{всей} системы соотношений. Они должны быть дополнены равенствами \eqref{eq18}, записываемыми как
\begin{align}
	\label{eq23}
	\begin{array}{ll}
		(\alpha d)_{k+q+\delta}=&(\beta c)_{k-q}\quad\forall k:\;0\leq k\leq n-p-1+q\\
		(ad)_{q+\delta-1}=&\sum\limits_{k=0}^{q+\delta-1}\epsilon^{q+\delta-1-k}(\alpha d)_k=0.
	\end{array}
\end{align}

Соотношения \eqref{eq19} имеют вид
\begin{align}
	\label{eq24}
	\begin{array}{ll}
		(\alpha c)_k<0 & \forall k:\;0\leq k\leq n-p-1\\
		(\beta d)_k<0 & \forall k:\;0\leq k\leq m-p-1
	\end{array}
\end{align}
и, фактически, выполняются в силу (\ref{eq21}, \ref{eq23}).

Наконец, вводя производящие функции
\begin{align}
	a(t)=\sum a_k(-t)^k\text{ и аналогичные для $b$, $c$ и $d$},
\end{align}
соотношения (\ref{eq21}, \ref{eq23}) могут быть записаны в виде
\begin{align}
	\label{eq26}
	(\alpha\beta cd)_{2k+\delta}>(2k+\delta)\left((\beta c)_{k-q}\right)^2\;\forall k: 0\leq 2k+\delta\leq n-p-1,
\end{align}
где $\alpha=(1-\epsilon t)a$, $\beta=(1-\epsilon t)b$, и
\begin{align}
	\label{eq27}
	ad=bc\cdot t^\Delta +Q,\;\deg Q\leq q+\delta-2\quad \text{для}\,\forall\deg(bc)\leq n-p-1.
\end{align}

Замечу, условие применения \eqref{eq21} несколько шире, чем \eqref{eq26}.


\section{Камни и корни}

В разделе доказано $S_n\mathop{\cap}S_m=\varnothing$ для $\forall\,m, n$.

Для краткости неполиномиальные части $R_n$ и $R_m$ функций $z_n$ и $z_m$ назовем их ядрами.

\noindent\textbf{Лемма (о полноте корней):} \textit{Любой корень уравнения $\lambda^{2p}=\Lambda_n$ представлен в ядре функции $z_n$ в ненулевом слагаемом $\mathrm{const}\cdot e^{i\lambda x}$.}

\noindent\textbf{Доказательство.} В силу четности и вещественности ядра вместе с каждым корнем $\lambda_j$ присутствует $-\lambda_j$ и $\overline{\mathstrut\lambda}_j$. Если в $R_n$ отсутствует $\lambda_0$, то из $L^{2p}(\sigma)$ в уравнении $L^{2p}R_n=0$ можно отбросить операторный множитель $\left(\sigma^2-\lambda_0^2\right)$, при отсутствии же $\lambda_j$ отбрасываем $\left(\sigma-\lambda_j\right)\left(\sigma-\overline{\lambda}_j\right)$; в итоге имеем $M(\sigma)R_n=0,\;\deg M\leq2p-2$; соответственно $L(\sigma)z_n=0,\;\deg L\leq2n-2$, где $L(\sigma)=M(\sigma)\cdot\sigma^{2n-2p}$. При отсутствии $\left(\sigma^2-\lambda_0^2\right)$ представим сам $L(\sigma)$, а иначе $\left(\sigma^2-\lambda_0^2\right)L(\sigma)$ в виде $\overline{\mathcal{A}}(\sigma)\mathcal{A}(\sigma)$, где $\mathcal{A}$ и $\overline{\mathcal{A}}$~--- комплексно-сопряженные полиномы от $\sigma$, вводимые аналогично \eqref{eq8}.

Имеем $\left\langle z_n\overline{\mathcal{A}} \mathcal{A} z_n\right\rangle=0,\;\deg \mathcal{A}=\deg\overline{\mathcal{A}}\leq n\Rightarrow z_n=0$ аналогично лемме о камне.

Пусть $\Lambda_n=\Lambda_m=\Lambda$. Представим:
\begin{align}
	R_n=\sum_{j=0}^{2p-1}r_{nj}e^{i\lambda_j x},\quad
	R_m=\sum_{j=0}^{2p-1}r_{mj}e^{i\lambda_j x};\quad
	\lambda_j^{2p}=\Lambda.
\end{align}
Обозначим $\mu=m-p-1$, $\nu=n-p-1$.

Сделав подстановку $x^2=1+(x^2-1)$ в полиномиальных частях $p_{2\nu}$ и $p_{2\mu}$, представим $z_n$ и $z_m$ в виде:
\begin{align}
	\label{eq30}
	z_n=R_n+\sum_{k=0}^\nu \gamma_{nk}(x^2-1)^k;\quad
	z_m=R_m+\sum_{k=0}^\mu \gamma_{mk}(x^2-1)^k.
\end{align}

Обозначим $\xi=x^2$, $\partial=\frac{d}{d\xi}$; если, как ранее, $d=\frac{d}{dx}$, то $d=2x\cdot\partial$. Индукцией по  $j$ легко доказывается:
\begin{align}
	\label{eq31}
	\left.z^{(i)}\right|_{x=1}=0\quad j=0, 1, \dots, l\Leftrightarrow\left.\partial\,^{j}z\right|_{\xi=1}=0\quad j=0, 1, \dots, l.
\end{align}

Из \eqref{eq30} получаем
\begin{align}
	\left.\partial\,^kz_l\right|_{\xi=1}=\left.\partial\,^kR_l\right|_{\xi=1} + \gamma_{lk}\cdot k!,\quad l=m, n.
\end{align}

Наконец, учитывая \eqref{eq31}, $z_n\in\Omega^n$, $z_m\in\Omega^m$, получим
\begin{align}
	\begin{array}{lll}
		\left.\partial\,^kR_n\right|_{\xi=1}=0 &\quad&
		k\in\mathcal{N}=\{n-p, \dots, n-1\} \\
		\left.\partial\,^kR_m\right|_{\xi=1}=0 &&
		k\in\mathcal{M}=\{m-p, \dots, m-1\}
	\end{array}
\end{align}


\section{Обсуждение}

В случае совпадения с.з. соотношения \eqref{eq18} или \eqref{eq27} задают довольно сильную связь между моментами и камнями функций $z_n$ и $z_m$. При наличии дополнительного исследования старших полиномиальных слагаемых в \eqref{eq5} или моментов \eqref{eq16} эти соотношения плюс неравенства \eqref{eq21} или \eqref{eq26}, либо прямо отношения (\ref{eq13}, \ref{eq15}) могут помочь как в решении вопроса о пересечении спектров для конкретных $\Delta=m-n$, так и проверить выполнимость гипотезы $S_n\mathop{\cap}S_m=\varnothing$ для любых $m, n$.

{\raggedleft\itshape
1--15 июля 2016 года

дер. Сондуга, г. Тотьма

}


\section{Благодарности}

Я выражаю свое восхищение А.\,И.\,Назарову, чей проницательный интерес к проблеме, хоть и опосредовано, был воспринят мной. Я очень признателен Юлии Петровой, которая передала для меня формулировку задачи в первоначальной постановке о минимуме $\Lambda_{n-1}<\Lambda_n$. Я благодарен Станиславу Крымскому за интересную дискуссию по работе и спасибо Павлу Муленко за согласие набрать данный текст.

Наконец, я выражаю свою сердечную благодарность Алексею, Ирине и Александру Завьяловым, Алексею Гущину и Любови Власовой за великолепные условия для работы и вдохновляющую поддержку.

\end{document}